\documentclass[a4paper,fleqn]{article}
\usepackage{authblk}
\usepackage[T1]{fontenc}
\usepackage[utf8]{inputenc}
\usepackage{lmodern} %
\usepackage[english]{babel}
\usepackage{hyperref}
\hypersetup{
	pdftitle={A new shape optimization approach for fracture propagation},
	pdfauthor={Tim Suchan, Kathrin Welker, and Winnifried Wollner}
}

\usepackage{color}
\usepackage{graphicx}
\graphicspath{{figures/},{svg-inkscape/},{tikzs/}}
\usepackage{subcaption}
\newlength\figureheight 
\newlength\figurewidth

\usepackage{amsmath,amssymb,amscd,latexsym,amsthm}
\usepackage{siunitx}
\usepackage{upgreek}
\usepackage{bbm} 
\usepackage{amsfonts}
\usepackage{mathrsfs}
\usepackage{cancel}
\usepackage{mathtools} %

\newcommand{\R}{{\mathbb{R}}} 

\newcommand{\toverbrace}[2]{\overbrace{#1}^{\textstyle#2}}

\DeclareMathOperator{\Div}{div}

\renewcommand{\vec}[1]{\boldsymbol{#1}}

\newcommand{\vn}{\vec{n}}
\newcommand{\state}{\vec{w}}
\newcommand{\shape}{u}
\newcommand{\fracture}{\rev{\mathfrak{U}}}
\newcommand{\stress}{\vec{\sigma}}
\newcommand{\strain}{\vec{\varepsilon}}
\newcommand{\hmetric}{\mathcal{H}}
\newcommand{\vV}{\vec{V}}

\newcommand{\Gb}{\varGamma_\text{bottom}}
\newcommand{\Gl}{\varGamma_\text{left}}
\newcommand{\Gr}{\varGamma_\text{right}}
\newcommand{\Gt}{\varGamma_\text{top}}

\newcommand{\dn}[1]{\frac{\partial #1}{\partial \vn}}

\newcommand{\rev}[1]{#1}

\begin{document}

\title{A new shape optimization approach for fracture propagation}

\author{Tim Suchan}
\affil{Helmut-Schmidt-Universität / Universität der Bundeswehr Hamburg, Holstenhofweg 85, 22043 Hamburg, Germany, \href{mailto:suchan@hsu-hh.de}{\ttfamily suchan@hsu-hh.de}\vspace*{6pt}}

\author{Kathrin Welker}
\affil{Technische Universität Bergakademie Freiberg, Akademiestraße 6, 09599 Freiberg, Germany, \href{mailto:Kathrin.Welker@math.tu-freiberg.de}{\ttfamily Kathrin.Welker@math.tu-freiberg.de}\vspace*{6pt}}

\author{Winnifried Wollner}
\affil{Universität Hamburg,  Bundesstr.~55, 20146 Hamburg, Germany, \href{mailto:winnifried.wollner@uni-hamburg.de}{\ttfamily winnifried.wollner@uni-hamburg.de}}

\date{}

\maketitle

\begin{abstract}
	Within this work, we present \rev{a novel approach} to fracture simulations
based on shape optimization techniques. Contrary to widely-used phase-field approaches
in literature the proposed method does not require a specified 'length-scale' parameter
defining the diffused interface region of the phase-field.
We provide the formulation and discuss the used solution approach.
We conclude with some numerical comparisons with well-established
single-edge notch \rev{tension and} shear tests.

\end{abstract}

\section{Introduction}
\label{sec:Introduction}
We consider an alternative approach for the solution of quasi-static
brittle fracture propagation due to Griffith's~\cite{Griffith:1921} model.
Based on the variational formulation proposed initially by Francfort \& Marigo~\cite{FrancfortMarigo:1998}.
In its simplest form this model consists of a minimization problem in displacement $\bar{\state}$ and
fracture $\bar{\fracture}$ which need to solve
\[
        \min_{\state,\fracture} E_{\rm{bulk}}(\state,\fracture) + G_c\mathcal{H}^{\rev{1}}(\fracture)
\]
over all admissible \rev{displacements $\state$} and fractures $\fracture$ which will be specified in
Section~\ref{sec:ProblemDescription}. Here $G_c \in \mathbb{R}$ denotes the fracture toughness and
$\mathcal{H}^{\rev{1}}$ is the $\rev{1}$-dimensional Hausdorff-measure. Due to the difficulty
of discretizing the lower-dimensional fracture a common approach is based on
Ambrosio \& Tortorelli~\cite{AmbrsosioTortorelli:1990}.
Here the lower-dimensional fracture is replaced by a, smooth,
phase-field function whose values indicate fractured or non fractured regions, see,
e.g.,~\cite{FrancfortMarigo:1998,BourdinFrancfortMarigot:2008,MieheWelschingerHofacker:2010,AmbatiGerasimovLorenzis:2015,NeitzelWickWollner:2018} for application to fracture problems.
The price to be paid in such phase-field problems is the introduction of (at least) two
regularization parameters, one for the approximation of the Hausdorff-measure and one to assert
the coercivity of the bulk-energy $E_{\rm{bulk}}$ in regions of vanishing phase-field.
To obtain meaningful numerical results the precise choice of the parameters and their balance with
the discretization error must be carefully addressed, see, e.g.,~\cite{WheelerWickWollner:2014},
and a necessary sharp resolution of the transition zone requires adaptive discretizations with
appropriate a posteriori error
indicators~\cite{BurkeOrtnerSuli:2010,Artina:2015,Mang2019,WallothWollner:2021}. 

Within this article, we propose a new approach avoiding the replacement of the lower dimensional
fracture by a phase-field method. The fracture evolution is then realized by means of techniques from shape optimization.
In order to obtain an efficient shape optimization algorithm, we consider the so-called Steklov-Poincar\'{e} metric \cite{SchulzSiebenbornWelker2015:2} in this work.
The Steklov-Poincar\'{e} metric has some numerical advantages over other types of metrics as shown in \cite{Siebenborn2017,Welker2016}.
In addition, the Steklov-Poincar\'{e} metric  allows to work with so-called weak formulations of shape derivatives, i.e., volume expression of shape derivatives. 
In the past, e.g., \cite{Delfour-Zolesio-2001,Sokolowski1991}, major effort in shape calculus has been devoted towards expressions for shape derivatives in the Hadamard form, i.e., in the boundary integral form. An equivalent and intermediate result in the process of deriving Hadamard expressions is a volume expression of the shape derivative, called the weak formulation. Thus, working with weak formulations saves analytical effort.

The rest of the paper is structured as follows, in Section~\ref{sec:ProblemDescription},
we will describe the considered problem in more detail.
In Section~\ref{sec:ShapeOptimizationApproach}, we will discuss how this problem can be restated
as a shape optimization problem and provide the formulas needed in the computation of
descent directions. Finally, in Section~\ref{sec:NumericalResults} we will provide numerical
results for the proposed approach for the well known single-edge notch tension and shear test
from~\cite{miehe2010phase} and the setup given in~\cite{Mang2019}.

\section{Problem description}
\label{sec:ProblemDescription}

\begin{figure}[tbp]
	\centering
	\setlength\figureheight{.5\textwidth} 
	\setlength\figurewidth{.5\textwidth}
	\newlength\svgwidth
	\setlength\svgwidth{\figurewidth}
	\centering%
	\begingroup%
  \makeatletter%
  \providecommand\color[2][]{%
    \errmessage{(Inkscape) Color is used for the text in Inkscape, but the package 'color.sty' is not loaded}%
    \renewcommand\color[2][]{}%
  }%
  \providecommand\transparent[1]{%
    \errmessage{(Inkscape) Transparency is used (non-zero) for the text in Inkscape, but the package 'transparent.sty' is not loaded}%
    \renewcommand\transparent[1]{}%
  }%
  \providecommand\rotatebox[2]{#2}%
  \newcommand*\fsize{\dimexpr\f@size pt\relax}%
  \newcommand*\lineheight[1]{\fontsize{\fsize}{#1\fsize}\selectfont}%
  \ifx\svgwidth\undefined%
    \setlength{\unitlength}{322.93000316bp}%
    \ifx\svgscale\undefined%
      \relax%
    \else%
      \setlength{\unitlength}{\unitlength * \real{\svgscale}}%
    \fi%
  \else%
    \setlength{\unitlength}{\svgwidth}%
  \fi%
  \global\let\svgwidth\undefined%
  \global\let\svgscale\undefined%
  \makeatother%
  \begin{picture}(1,1.03556101)%
    \lineheight{1}%
    \setlength\tabcolsep{0pt}%
    \put(0,0){\includegraphics[width=\unitlength,page=1]{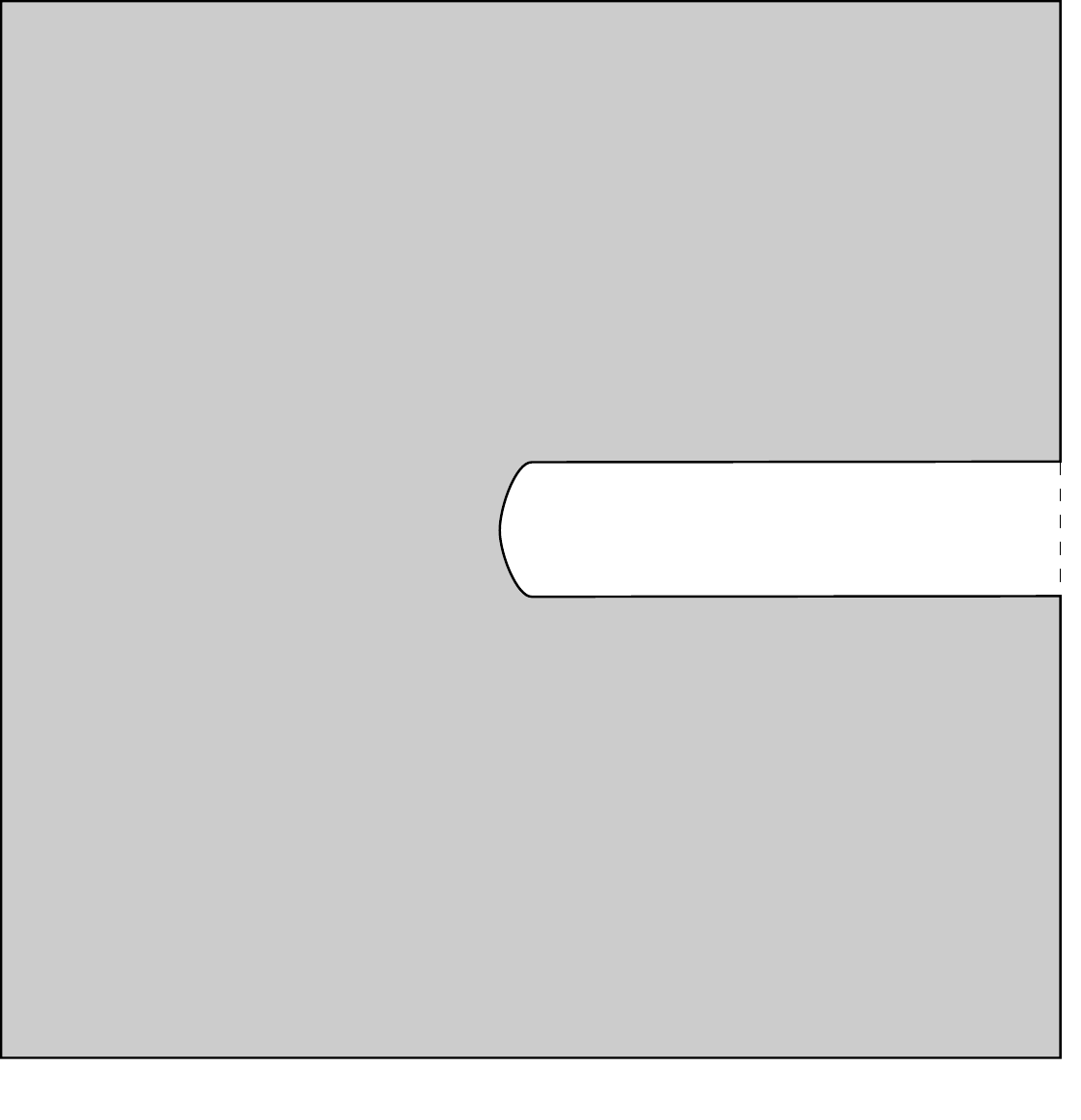}}%
    \put(0.93286519,0.64646842){\color[rgb]{0,0,0}\makebox(0,0)[t]{\lineheight{1.25}\smash{\begin{tabular}[t]{c}$\color{blue}P_1$\end{tabular}}}}%
    \put(0.54221015,0.62671695){\color[rgb]{0,0,0}\makebox(0,0)[t]{\lineheight{1.25}\smash{\begin{tabular}[t]{c}$\color{red}u$\end{tabular}}}}%
    \put(0.92654376,0.42010833){\color[rgb]{0,0,0}\makebox(0,0)[t]{\lineheight{1.25}\smash{\begin{tabular}[t]{c}$\color{blue}P_2$\end{tabular}}}}%
    \put(0.72057128,0.51895167){\color[rgb]{0,0,0}\makebox(0,0)[t]{\lineheight{1.25}\smash{\begin{tabular}[t]{c}$\widetilde{\varOmega}$\end{tabular}}}}%
    \put(0.24101896,0.85512873){\color[rgb]{0,0,0}\makebox(0,0)[t]{\lineheight{1.25}\smash{\begin{tabular}[t]{c}$\varOmega$\end{tabular}}}}%
    \put(0,0){\includegraphics[width=\unitlength,page=2]{initial_crack_svg-tex.pdf}}%
  \end{picture}%
\endgroup%

	\caption{Sketch of the hold-all domain $D$.}
	\label{fig:sketchModelDescription}
	\vspace*{-2mm}
\end{figure}

We focus on a two-dimensional setup. Here a hold-all domain $D \subset \R^2$
is considered which is decomposed into a fracture $\fracture$ and a remaining domain $\varOmega\subset \mathbb{R}^2$ such that $D=\varOmega \sqcup \fracture$. A sketch can be found in Fig.~\ref{fig:sketchModelDescription}. 

Given an initial fracture $\fracture^0$ the quasi-static evolution
of the fracture $\fracture$ is governed by the energy minimization 
\begin{align}
	\min_{\state,\fracture} \frac{1}{2} \left( \mathbb{C} : \strain(\state), \strain(\state) \right)_{L^2(\varOmega)} - (\vec{f},\state)_{L^2(\varOmega)} + G_c \hmetric^{\rev{1}}(\fracture).
	\label{eqn:MinimizationWithRespectToStateAndShape}
\end{align}
Here $\state \in H^1_D(\varOmega;\R^{\rev{2}}) + \state^D = \{\state \in H^1(\varOmega;\R^{\rev{2}}) \,|\, \state = \state^D \text{ on }\varGamma_D\}$ are the admissible displacements, where $\varGamma_D \subset \partial \varOmega \setminus \fracture$ is a given boundary part, and $\state^D$ denotes prescribed Dirichlet data.
\rev{The fracture $\fracture$ is required to be monotonically increasing in time, i.e., $\fracture(t+\delta t) \supset \fracture(t)$ for any $\delta t > 0$.}
The first two terms in~\eqref{eqn:MinimizationWithRespectToStateAndShape} denote the elastic
bulk energy, i.e.,
$\strain(\state)=\frac{1}{2} \left( \nabla \state + \left( \nabla \state \right)^\top \right)$
is the symmetric gradient, and the stress-strain relation is given by
$\stress(\state) = \mathbb{C} : \strain(\state) = 2 \mu \strain(\state) + \lambda \operatorname{tr}(\strain(\state)) I$ with the Lam\'e parameters $\mu,\lambda$.
The last term  in~\eqref{eqn:MinimizationWithRespectToStateAndShape} is the $\rev{1}$-dimensional Hausdorff measure which describes the length of the fracture. 

In order to solve (\ref{eqn:MinimizationWithRespectToStateAndShape})  with shape optimization techniques, we need to find a shape definition of the fracture. One possibility is to  define the fracture by a curve
$\shape\colon [0,1]\to \mathbb{R}^2$ 
and assume that $\shape$ is an element of
\begin{equation*}
	\label{B_e_0}
	B_e^0=B_e^0\left([0,1],\R^2\right)\coloneqq\mbox{Emb}^0([0,1],\R^2)/\mbox{Diff}^0([0,1]),
\end{equation*}
where
\begin{equation*}
	\begin{split}
		\mbox{Emb}^0([0,1],\R^2) \coloneqq \{&\phi\in C^\infty([0,1],\R^2)\colon
		\phi(0)=P_1,\ \phi(1)=P_2,\,\\
		&\phi \mbox{ proper injective immersion}\}, \quad \text{and} \\
		\mbox{Diff}^0([0,1]) \coloneqq \{&\phi\colon [0,1]\to[0,1]\colon 
		\phi(0)=0,\ \phi(1)=1, \ \phi \mbox{ diffeomorphism}\}.
	\end{split}
\end{equation*}
\rev{Here, $P_1$ and $P_2$ are the prescribed start and end of the curve.}
The shape space  $B_e^0\left([0,1],\R^2\right)$ is constructed in analogy to the manifold $B_e(S^1,\R^2)$ introduced by Michor \& Mumford in \cite{Michor2006}. 
In this paper, we focus on such a shape definition because the shape
space $B_e(S^1,\R^2)$ is well-investigated (cf., e.g.,
\cite{KrieglMichor}) and commonly used in the connection with shape
optimization in the last years
\cite{GeiersbachLoayzaWelker,Siebenborn2017,Schulz,LuftWelker,SchulzWelker,SchulzSiebenbornWelker2015:2,Welker_diffeological}.
\rev{This induces a new splitting of the domain $D$ into the
physically remaining domain $\varOmega\subset \mathbb{R}^2$, the fracture boundary curve
$\shape$ and a remainder domain $\widetilde{\varOmega}\subset \mathbb{R}^2$ as $D=\varOmega \sqcup \fracture \sqcup \widetilde{\varOmega}$. A sketch can be found in Fig.~\ref{fig:sketchModelDescription}. Corresponding to the monotonically increasing fracture, we now require the domain $\widetilde{\varOmega}$ to be monotonically increasing in time, i.e., $\widetilde{\varOmega}(t+\delta t) \supset \widetilde{\varOmega}(t)$ for any $\delta t > 0$.}

The $\rev{1}$-dimensional Hausdorff measure of the
fracture $\fracture$ is approximated by $\frac{1}{2}
\hmetric^{\rev{1}}(\shape)$, since $\hmetric^{\rev{1}}(\shape)$ describes the
length of the curve $\shape$ which covers the two boundaries of the
fracture \rev{and is therefore twice as long}. Additionally, it has to be incorporated that the \rev{domain $\widetilde{\varOmega}$ may not increase} in volume. To account for this, a mild volume regularization
$J_{reg}(\shape) = \nu \cdot \int_{\widetilde{\varOmega}} 1 \, \mathrm{d} \vec{x}$ is added to the objective functional.
Here, $\nu>0$ is a parameter which is chosen s.t. the term has a minuscule effect in comparison to the original objective functional.

Equation~\eqref{eqn:MinimizationWithRespectToStateAndShape} can be reformulated to a minimization w.r.t. to the shape~$\shape$ subject to a partial differential equation (PDE) constraint:
\vspace*{-.3cm}
\begin{align}
	\min_{\shape\in B_e^0} \toverbrace{
		\toverbrace{
			\frac{1}{2} \left( \vec{\sigma}(\state), \vec{\varepsilon}(\state) \right)_{L^2(\varOmega)} -(\vec{f},\state)_{L^2(\varOmega)}
		}{E_{\text{bulk}}} + \toverbrace{
			\frac{1}{2} G_c \hmetric^{\rev{1}}(\shape)
		}{E_{\text{fracture}}} + \toverbrace{\nu \cdot \int_{\widetilde{\varOmega}} 1 \, \mathrm{d} \vec{x}}{J_{reg}}
	}{=J(\shape)}
\label{eqn:MinimizationWithRespectToShape}
\end{align}
\vspace*{-.5cm}
\begin{subequations}
	\begin{align}
		\text{ s.t. } - \Div(\stress(\state)) &= \vec{f} &&\text{in } \varOmega \\
		\state &= \state^D &&\text{on } \varGamma^D&&&&& \\
		\vec{\sigma}(\state) \, \vec{n} &= \vec{0} &&\text{on } \varGamma^N = \partial\varOmega \setminus \varGamma^D.&&&&&
	\end{align}
	\label{eqn:PDEConstraintGeneralForm}
\end{subequations}
Here, commonly-used boundary conditions are already included. For further details, cf. e.g.~\cite{AmbatiGerasimovLorenzis:2015,Mang2019}.

\section{Shape optimization approach}
\label{sec:ShapeOptimizationApproach}

In classical shape calculus, a shape is considered to be a subset of $\mathbb{R}^{\rev{2}}$, only. However, equipping a shape with additional structure provides theoretical advantages, enabling the use of concepts from differential geometry. In our setting, a shape is assumed to be an element of the shape manifold $B_e^0$. 
Because of the equivalence relation $\mbox{Diff}^0([0,1])$, the tangent space is isomorphic to the set of all smooth vector fields along $u$, i.e.,
$
T_uB_e^0\left([0,1],\R^2\right)\cong\{\vec{h}\colon \vec{h}=\alpha \vn, \alpha\in C^\infty \hspace{-.4mm}\left([0,1]\right)\},
$
where $\vn$ is the unit outward normal of $\varOmega$ at $u$. 

In view of obtaining gradient-based optimization approaches, the gradient needs to be specified. 
In addition, one needs to use a gradient-based optimization approach in Riemannian manifolds which can be found for example in \cite{GeiersbachHandbook,Welker_diffeological}.
The gradient will be characterized by the chosen Riemannian metric on $B_e^0$. As already mentioned above, we focus on the Steklov-Poincar\'{e} metric in this paper.
The detailed definition of the shape gradient with respect to the metric can be found in \cite[Definition 4]{GeiersbachHandbook}. 
Roughly spoken, in order to compute the shape gradient with respect to the Steklov-Poincar\'{e} metric, one needs to solve the so-called deformation equation
$
a(\vec{V},\vec{W})=\mathrm{d}J(\shape)[\vec{W}] $ for all $ \vec{W}\in H_V^1(\varOmega,\R^2)\cup C^\infty(\varOmega,\R^2),
$
where $a(\cdot,\cdot)$ denotes a coercive and symmetric bilinear form, $\mathrm{d}J(\shape)[\vec{W}]$ is the shape derivative of $J(\shape)$ in the direction $\vec{W}$, and $H^1_V(\varOmega, \R^{\rev{2}})=\left\{ \vV \in H^1(\varOmega, \R^{\rev{2}}) \,|\, \vV = \vec{0} \text{ on } \partial \varOmega \setminus u \right\}$.
In our numerics, the bilinear form is chosen to be the weak form of the linear elasticity equation.
Following standard techniques for the calculation of shape derivatives, we get
\begin{align}
\begin{aligned}
	\mathrm{d}J(\shape)[\vec{W}] = &\int_{\varOmega(\shape)} - \frac{1}{2} \left( \nabla \state \nabla \vec{W} + \left( \nabla \state \nabla \vec{W} \right)^\top \right) : \stress(\state) - \vec{W}^\top \nabla \vec{f} \, \state \\
	&\hphantom{\int_{\varOmega(\shape)}\,} + \Div(\vec{W}) \cdot \left( \frac{1}{2} \stress(\state) : \strain(\state) - \vec{f}^\top \state \right) \mathrm{d} \vec{x} \\
	&+ \frac{1}{2} G_c \int_{\shape} \kappa \vec{W}^\top \vn \,\mathrm{d} s - \int_\shape \nu \vec{W}^\top \vn \,\mathrm{d} s,
\end{aligned}
\end{align}
where $\kappa$ denotes the curvature of the respective curve. %
The resulting $\vV$, the so-called mesh deformation, is then modified to avoid shrinking of the fracture even when applying low loads, as already described in section~\ref{sec:ProblemDescription}. For this, we require $\left< \vV\vert_u, \vn \right> \leq 0$ (pointwise) nearly everywhere on $u$. However, since the mesh deformation is defined in all of $\varOmega$, an extension of the normal vector is  beneficial from a numerical point of view. The extension is performed by solving the Eikonal equation
\begin{align}
	\label{eqn:EikonalEquation}
\begin{aligned}
	\left| \nabla \Phi(\vec{x}) \right| = 1 \quad \text{in } \varOmega , \quad
	\Phi(\vec{x}) = 0 \quad \text{on } u , \quad
	\dn{\Phi(\vec{x})} = 0 \quad \text{on } \partial \varOmega \setminus u.
\end{aligned}
\end{align}
The iterative solution process is stabilized by adding a linear second-order derivative term to equation~\eqref{eqn:EikonalEquation}. %
The extended normal vector into the domain can then be computed by $\vec{N} = \nabla \Phi(\vec{x})$. Using this extended normal vector, we require $\left< \vV, \vec{N} \right> \leq 0$ (pointwise) nearly everywhere in $\varOmega$.

\section{Numerical results}
\label{sec:NumericalResults}

We study the fracture propagation for two common benchmark problems, cf., e.g.,~\cite{AmbatiGerasimovLorenzis:2015,Mang2019}: 
 the single-edge notched tension test using linear elasticity with a homogeneous material and
 the single-edge notched shear test, also using linear elasticity with a homogeneous material.
The material is described by the Lamé parameters $\lambda = \SI{121.15e3}{\newton\per\square\milli\meter}$ and $\mu = \SI{80.77e3}{\newton\per\square\milli\meter}$. \rev{To ensure that $\left< \vV, \vec{N} \right> \leq 0$ is fulfilled pointwise nearly everywhere in $\varOmega$ in the numerical computations, we compute the mesh deformation~$\vV$ as described in section~\ref{sec:ShapeOptimizationApproach} without constraints, but then replace $\vV$ at any point where the constraint is violated by $\vV=\vec{0}$.} The values of bulk energy $E_{\text{bulk}}$ as shown in equation~\eqref{eqn:MinimizationWithRespectToShape} and of the boundary force $\vec{\tau} = \int_{\Gt} \stress(\state) \, \vn \,\mathrm{d} s$ versus the applied displacement $\state^D$ on $\Gt$ are used to compare the results to previously-published literature.

\begin{figure}[tb]
	\centering%
	\setlength\figureheight{.45\textwidth} %
	\setlength\figurewidth{.45\textwidth}%
	\newlength\svgwidth
	\setlength\svgwidth{\figurewidth}
	\begin{subfigure}[b]{\figurewidth}%
		\centering%
		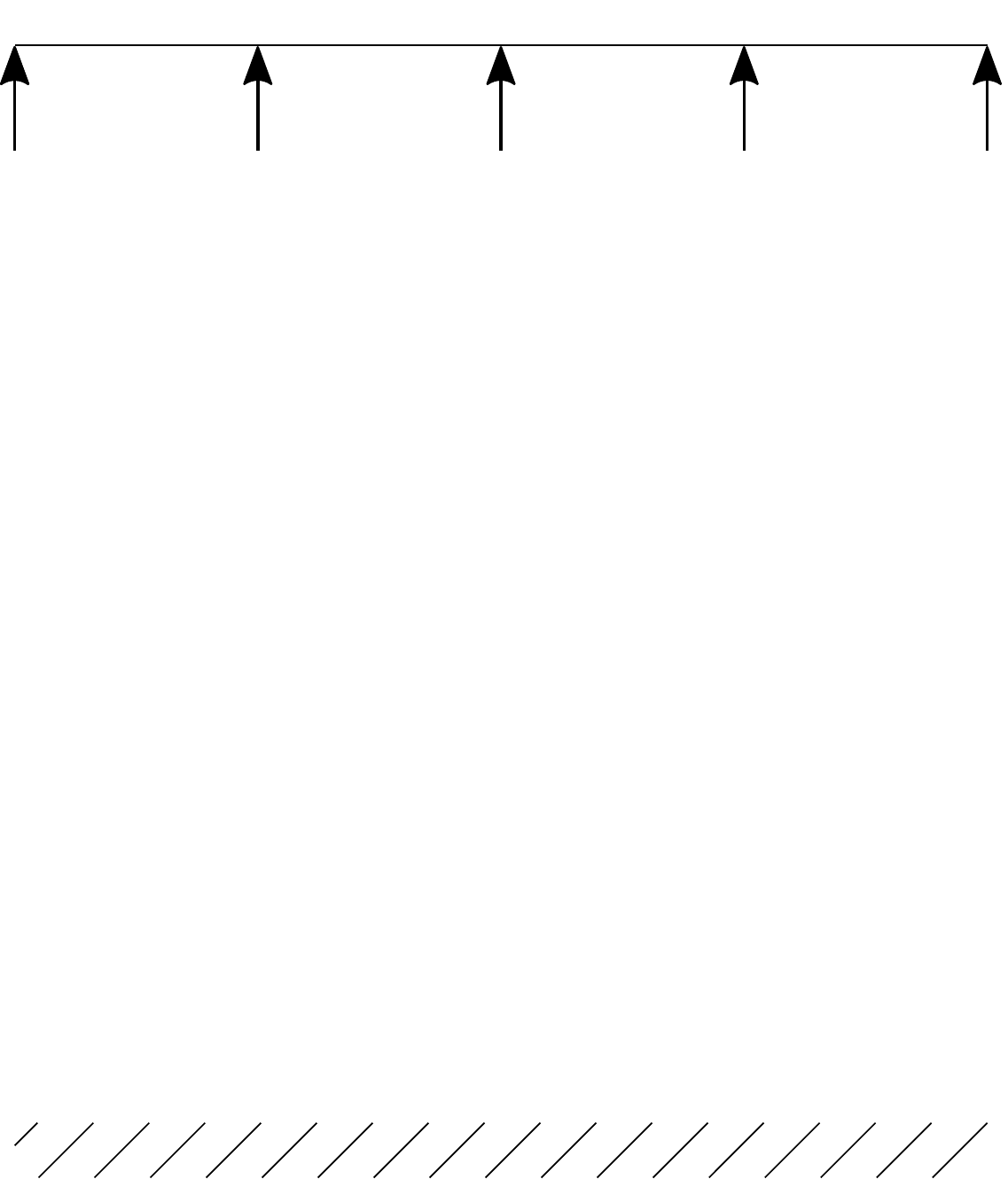
	\end{subfigure}%
	\quad%
	\newlength\svgwidth
	\setlength\svgwidth{\figurewidth}
	\begin{subfigure}[b]{\figurewidth}%
		\centering%
		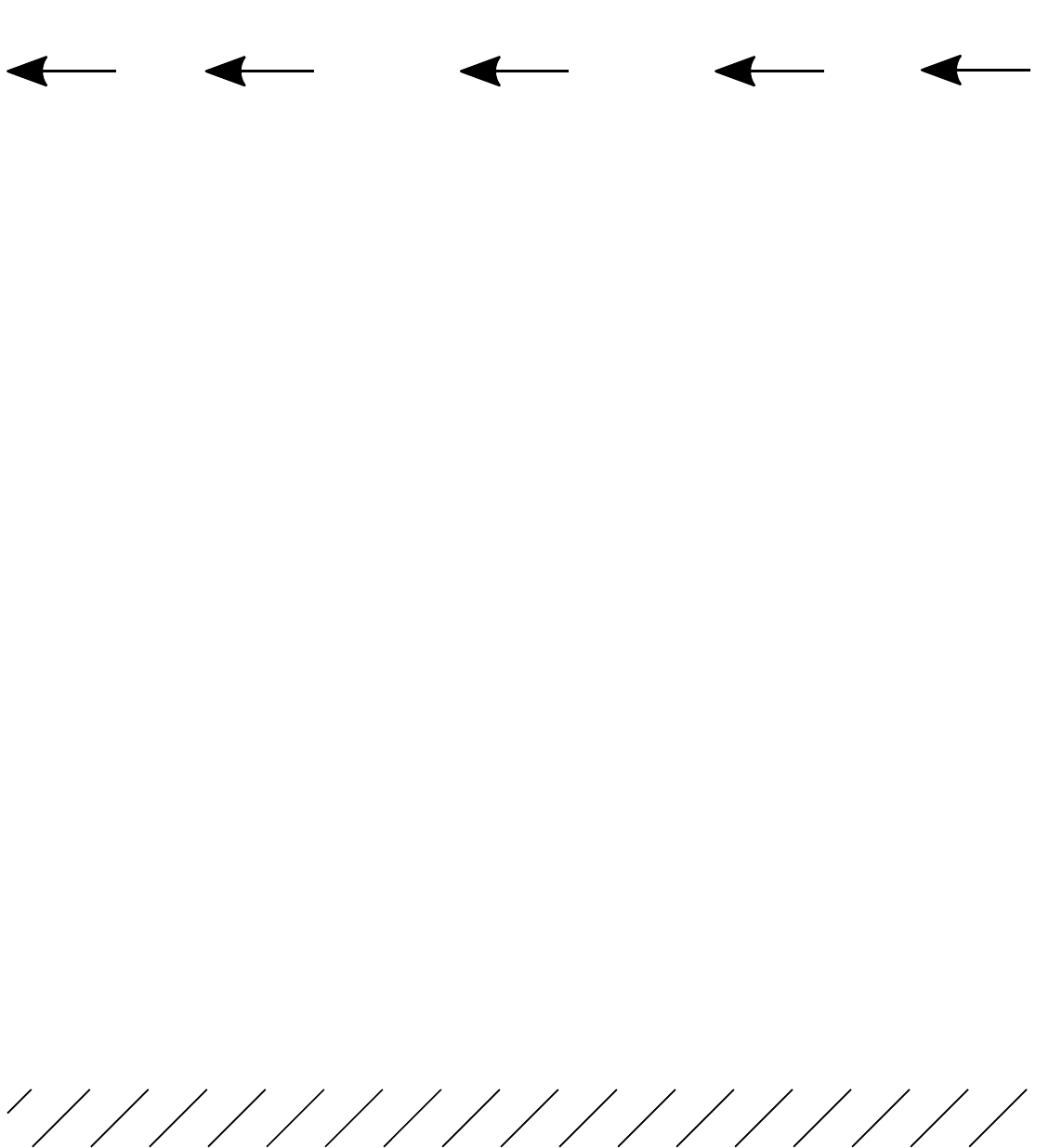
	\end{subfigure}%
	\caption{Sketches of the computational domain~$\varOmega$ with boundaries and boundary conditions for the tension test (left) and for the shear test (right).}
	\label{fig:sketchModelDescriptionWithBoundaries}%
	\vspace*{-2mm}
\end{figure}

Unless otherwise described, the length measures are given in $\si{\milli\meter}$. For the numerical experiments, the hold-all domain $D$ is chosen as $D=\left( 0, 1 \right)^2$, with the start point of the curve $\shape^0$ at $P_1=\left( 1, 0.51 \right)^\top$, extending to $\left( 0.5, 0.51 \right)^\top$, followed by a rounded tip with a radius of $0.01$, and continuing from $\left( 0.5, 0.49 \right)^\top$ until $P_2=\left( 1, 0.49 \right)^\top$. The domain $\widetilde{\varOmega}$ is not required for the numerical computations, thus we can restrict our numerical investigations on the computational domain~$\varOmega$, cf. Fig.~\ref{fig:sketchModelDescriptionWithBoundaries}.
The boundary of computational domain is split into separate parts as shown in Fig.~\ref{fig:sketchModelDescriptionWithBoundaries}. These parts are assigned into the two disjoint sets $\varGamma_D=\Gb \cup \Gt$ and $\varGamma_N= \Gl \cup \Gr \cup u$. In both, the tension and the shear test, a homogeneous Dirichlet boundary condition on $\Gb$, i.e. $\state^D = \vec{0} \quad \text{on } \Gb$, and a homogeneous Neumann boundary condition on $\varGamma_N$ are imposed. Since volumetric forces are neglected in the benchmark problems, the volumetric force term is set to $\vec{f}=\vec{0}$.

The discretization of the domain is performed using Gmsh \cite{Geuzaine2009}, and yields 1220 nodes and 2438 triangular elements initially. Linear Finite Elements are used throughout this work. Due to the large deformations of the shape an automatic remesher is activated when the mesh quality deteriorates past a threshold. A finer discretization is used near the fracture tip to resolve the stress singularity with better accuracy. Because of the remeshing, the number of nodes and elements changes throughout the optimization.

For the minimization of equation~\eqref{eqn:MinimizationWithRespectToShape} subject to the PDE constraint~\eqref{eqn:PDEConstraintGeneralForm}, a gradient descent algorithm with a constant step size of $10^{-2}$ is implemented in FEniCS 2019.1.0 \cite{Alnaes2015}. The stop criterion for the gradient descent is fulfilled if the change in length of the curve~$u$ is below a value of $10^{-8}$.

\paragraph{Tension test.}

For the tension test, the volume regularization parameter is chosen as $\nu=1.0$. 
The prescribed Dirichlet data on $\Gt$ for the tension test is chosen as $\state^D = t \cdot \left( 0, 10^{-5} \right)^\top$ as shown in Fig.~\ref{fig:sketchModelDescriptionWithBoundaries} on the left, where $t=1,2,\ldots$ is increased until fracture. To reduce the computational effort while the fracture has not started propagating, the displacement increment is set to $10^{-3}$ until a displacement of $4 \cdot 10^{-3}$ is reached.

\begin{figure}[tb]
	\centering
	\setlength\figureheight{.3\textheight} 
	\setlength\figurewidth{.3\textwidth}
	\begin{subfigure}[t]{\figurewidth}
		\centering
		\includegraphics[width=\figurewidth]{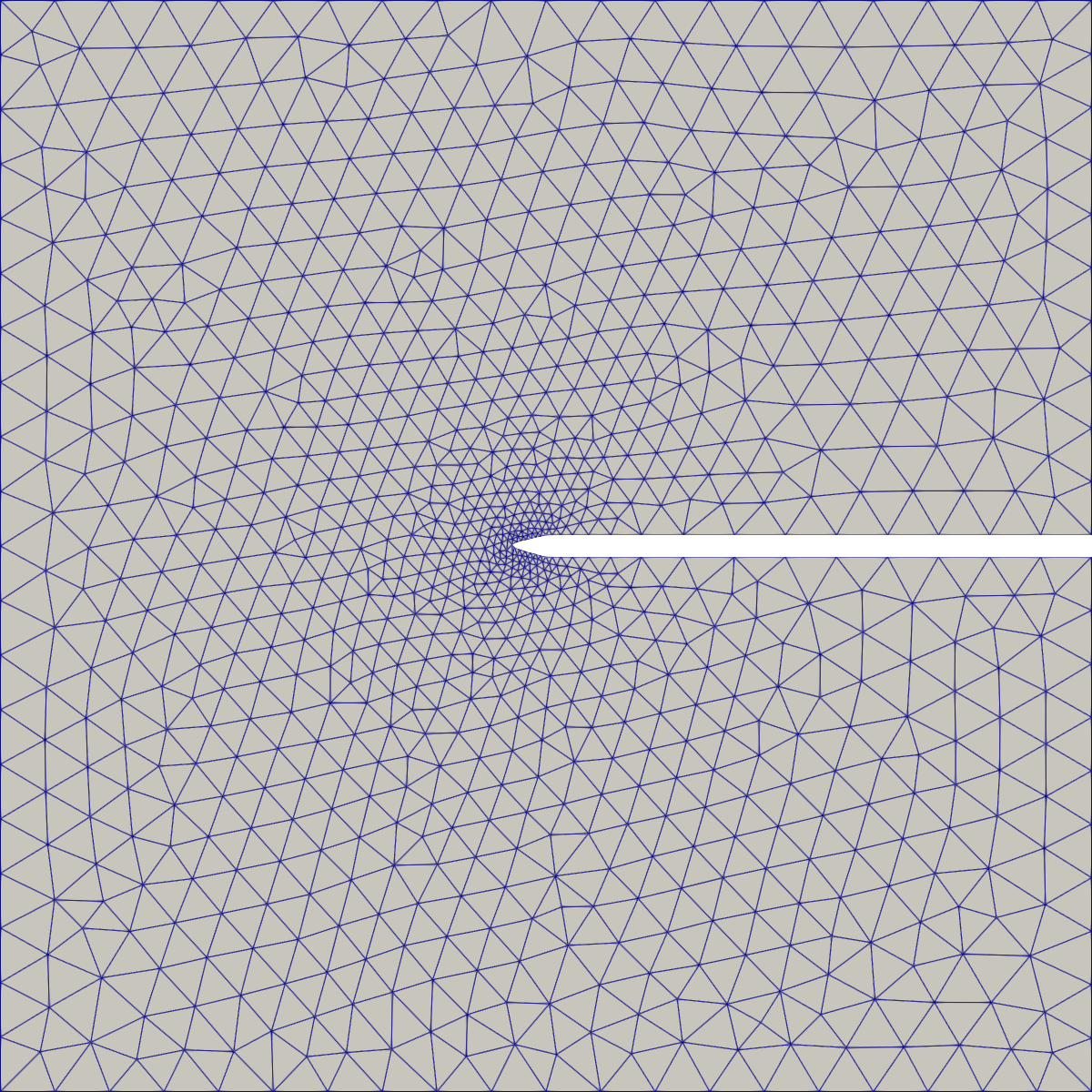}
		\caption{Starting from an initial rounded tip, a more sharp tip has been forming, here shown at $\SI{4.44e-3}{\milli\meter}$.}
		\label{fig:mesh_tension_initial}
	\end{subfigure}
	\quad
	\begin{subfigure}[t]{\figurewidth}
		\centering
		\includegraphics[width=\figurewidth]{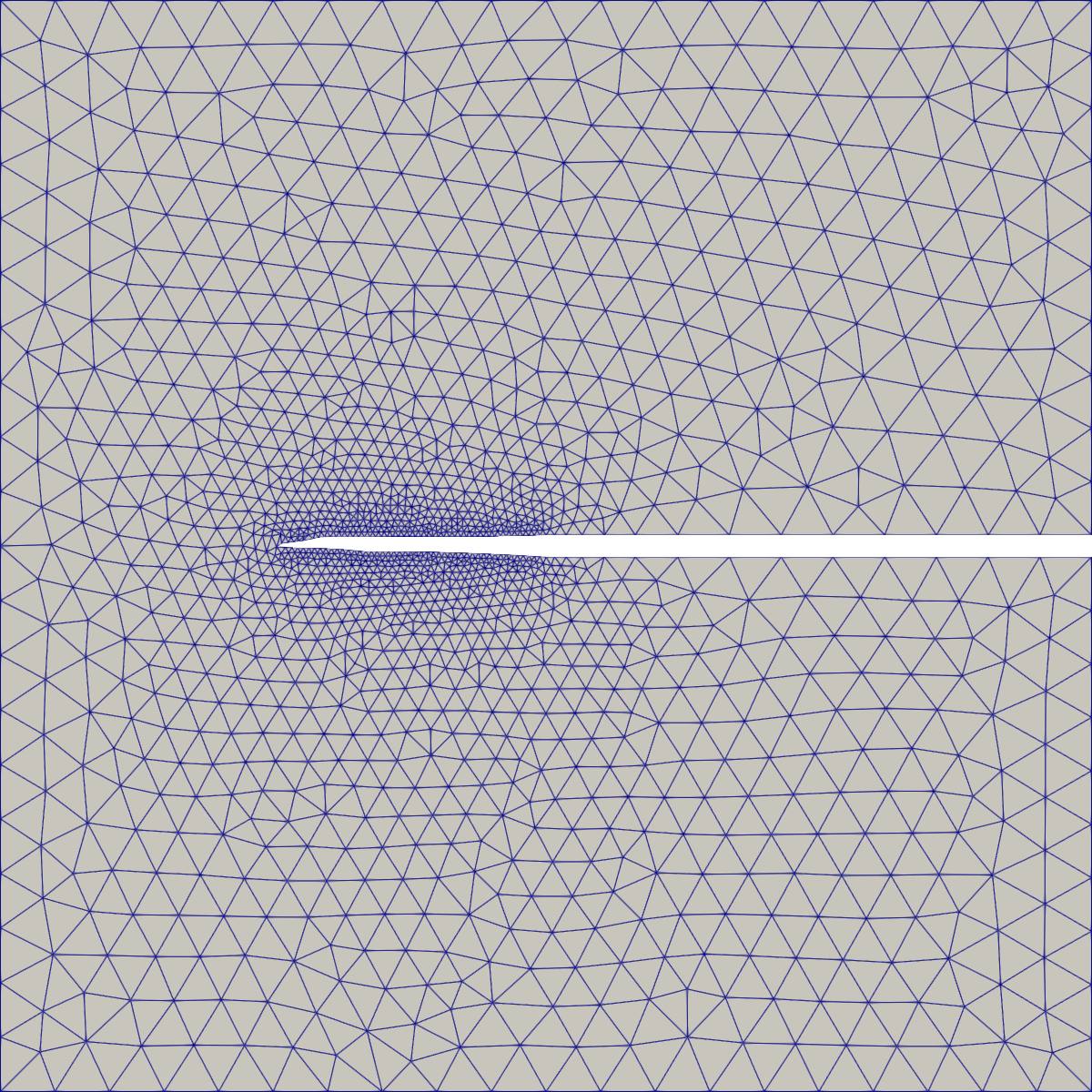}
		\caption{The fracture has already started propagating at $\SI{4.67e-3}{\milli\meter}$.}
		\label{fig:mesh_tension_middle}
	\end{subfigure}
	\quad
	\begin{subfigure}[t]{\figurewidth}
		\centering
		\includegraphics[width=\figurewidth]{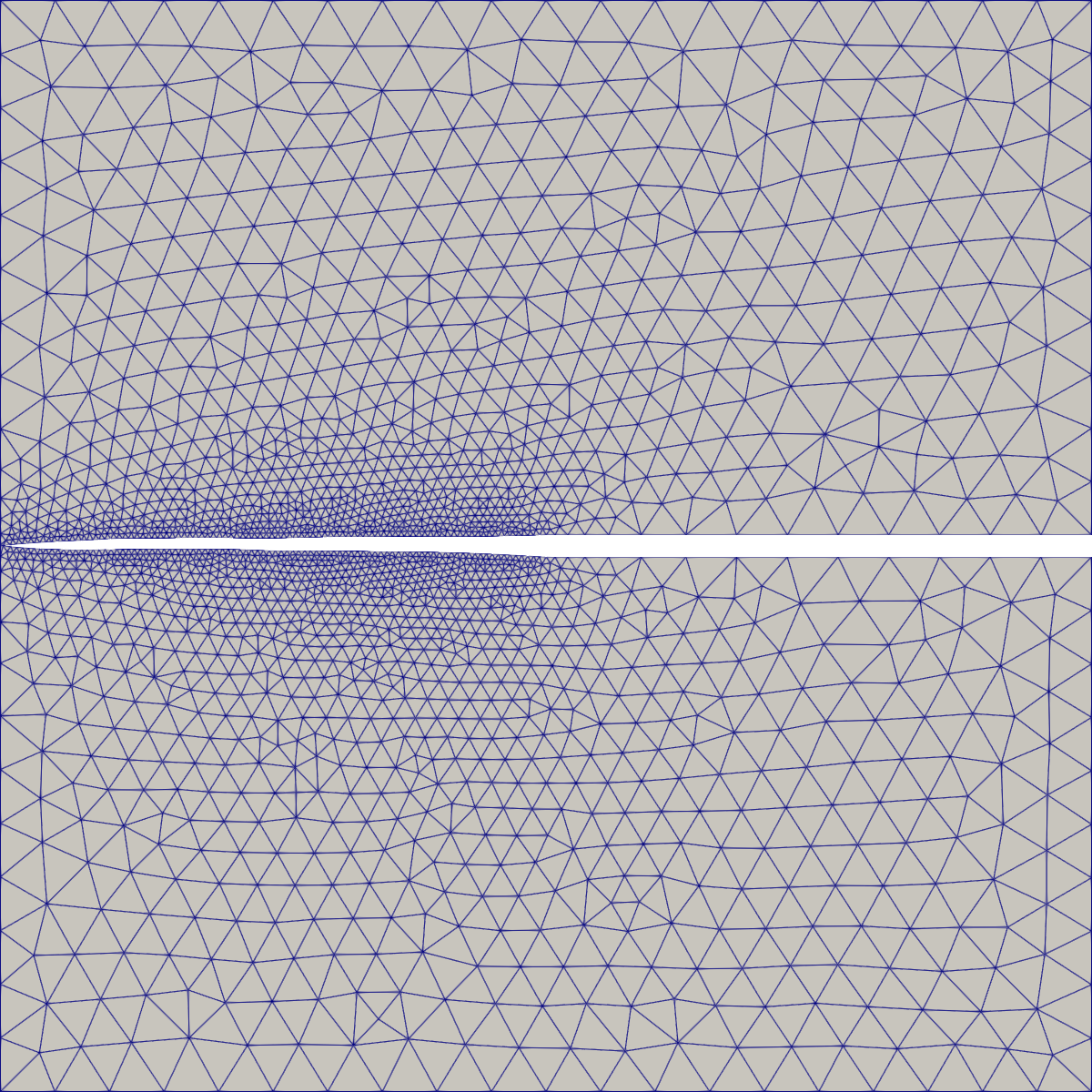}
		\caption{The fracture completely propagated through the domain at $\SI{4.87e-3}{\milli\meter}$.}
		\label{fig:mesh_tension_end}
	\end{subfigure}
	\caption{Fracture propagation for the tension test.}
	\label{fig:mesh_tension}
	\vspace*{\floatsep}%
	\centering%
	\setlength\figureheight{5cm}%
	\setlength\figurewidth{.45\textwidth}%
	\begin{subfigure}[t]{\figurewidth}%
		\centering%
		\includegraphics{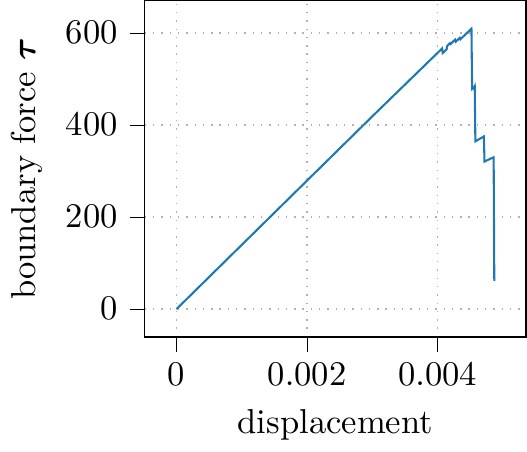}
		\label{fig:tension_boundary_force}%
	\end{subfigure}%
	\quad%
	\begin{subfigure}[t]{\figurewidth}%
		\centering%
		\includegraphics{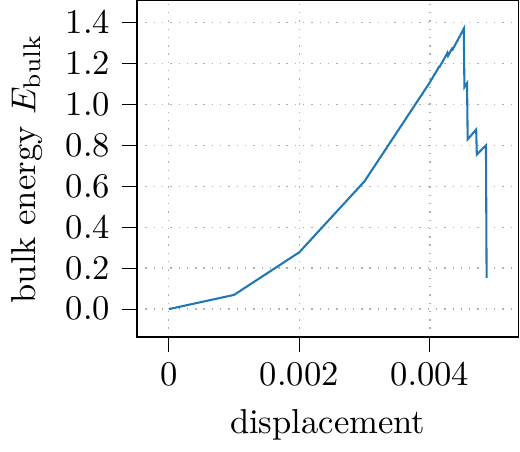}
		\label{fig:tension_bulkEnergy}%
	\end{subfigure}%
	\caption{Plot of boundary force $\vec{\tau}$ (left) and bulk energy $E_{\text{bulk}}$ (right) over displacement for the tension test.}%
	\label{fig:tension}%
	\vspace*{-2mm}
\end{figure}

Fig.~\ref{fig:mesh_tension} shows the fracture for different applied displacements. At a displacement of $4\cdot 10^{-3}$, the originally round tip of the fracture starts to become pointier, gradually becoming slightly longer with further increasing displacements. An example for this is shown in Fig.~\ref{fig:mesh_tension_initial} for a displacement of $4.44\cdot 10^{-3}$. The first major fracture growth occurs in a horizontal direction at $4.53\cdot 10^{-3}$. The propagating fracture at $4.67 \cdot 10^{-3}$ is shown in Fig.~\ref{fig:mesh_tension_middle}. At an applied displacement of $4.87 \cdot 10^{-3}$, the fracture has propagated horizontally through the whole computational domain, as illustrated in Fig.~\ref{fig:mesh_tension_end}. The mesh in the last iteration contains 2076 nodes and 4150 elements.

Fig.~\ref{fig:tension} illustrates the boundary force $\vec{\tau}$ (left) and bulk energy $E_{\text{bulk}}$ (right) over the applied Dirichlet data. Until a displacement of $4 \cdot 10^{-3}$, a linear increase in boundary force and a quadratic increase in bulk energy can be observed. In this regime, no fracture growth has happened yet. From $4 \cdot 10^{-3}$ until $4.5 \cdot 10^{-3}$, the boundary force is still increasing, however the slight growth of the fracture tip is causing a \rev{slightly slower increase} in boundary force and bulk energy than before. At $4.53\cdot 10^{-3}$, the fracture then rapidly grows in length, which can be spotted in both boundary force and bulk energy. The fracture has propagated fully through the computational domain at $4.87 \cdot 10^{-3}$.

These results are in good accordance with the previously-mentioned literature, however we have detected a fracture propagation at slightly lower displacements than the literature, which is most likely caused by the different algorithm. The results in the literature are generated with only few iterations per displacement, whereas for the results for this publication we have performed a gradient descent for each displacement until the stop criterion mentioned in the introduction to this section.

\paragraph{Shear test.}

For the shear test, the volume regularization parameter is set to $\nu=10.0$ and $\state^D = t \cdot \left( -10^{-5} , 0 \right)^\top$, which is illustrated in Fig.~\ref{fig:sketchModelDescriptionWithBoundaries} on the right, with increasing $t=1,2,\ldots$ until fracture. Similar to the tension test, the displacement increment is set to $10^{-4}$ until a displacement of $8 \cdot 10^{-3}$ is reached to reduce the computational effort.

\begin{figure}[tb]
	\centering%
	\setlength\figureheight{.3\textheight}%
	\setlength\figurewidth{.3\textwidth}%
	\begin{subfigure}[t]{\figurewidth}%
		\centering%
		\includegraphics[width=\figurewidth]{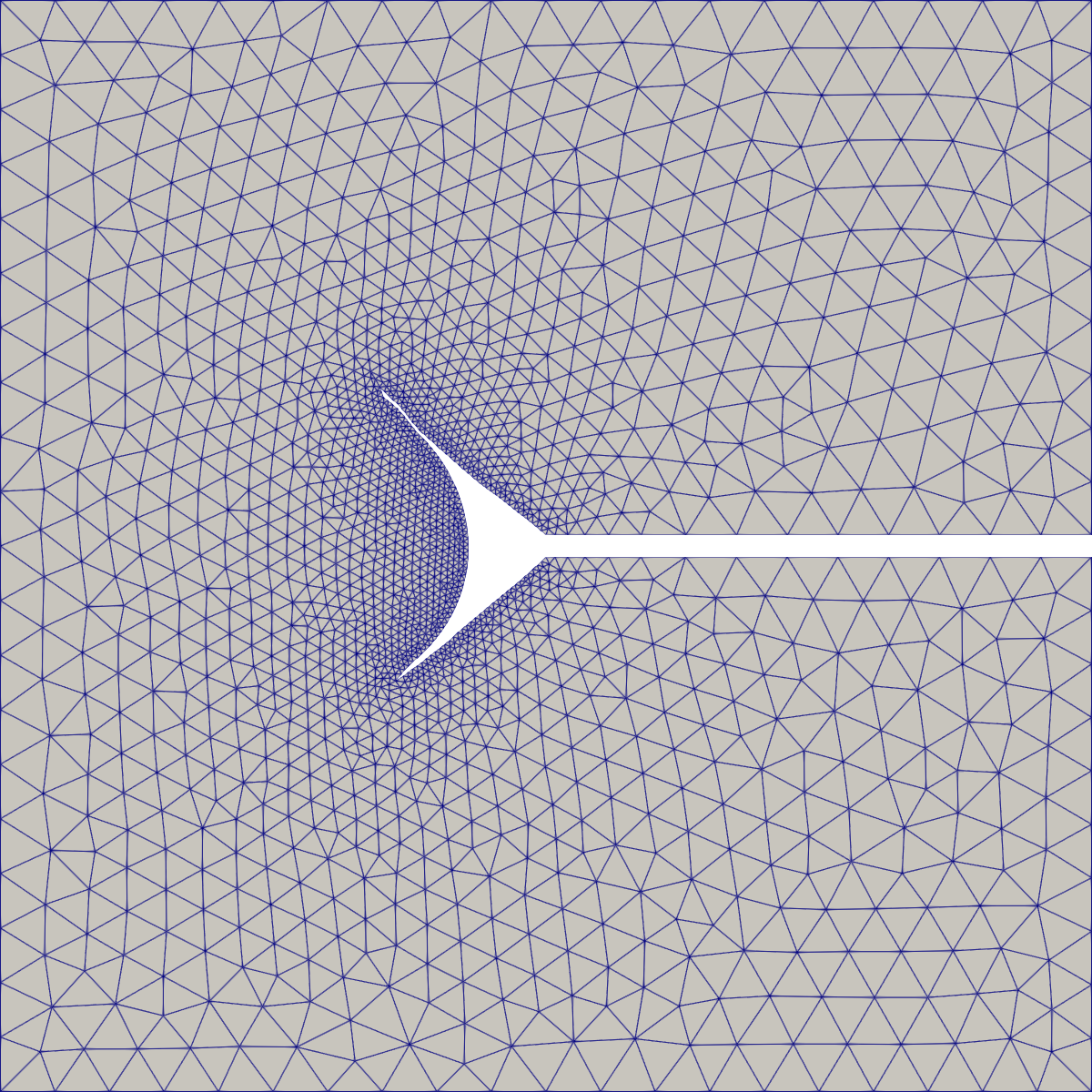}%
		\caption{The fracture at the end of the gradient descent for a displacement of $\SI{9.26e-3}{\milli\meter}$, the first major propagation of the fracture, is shown here.}%
		\label{fig:mesh_shear_initial}%
	\end{subfigure}%
	\quad%
	\begin{subfigure}[t]{\figurewidth}%
		\centering%
		\includegraphics[width=\figurewidth]{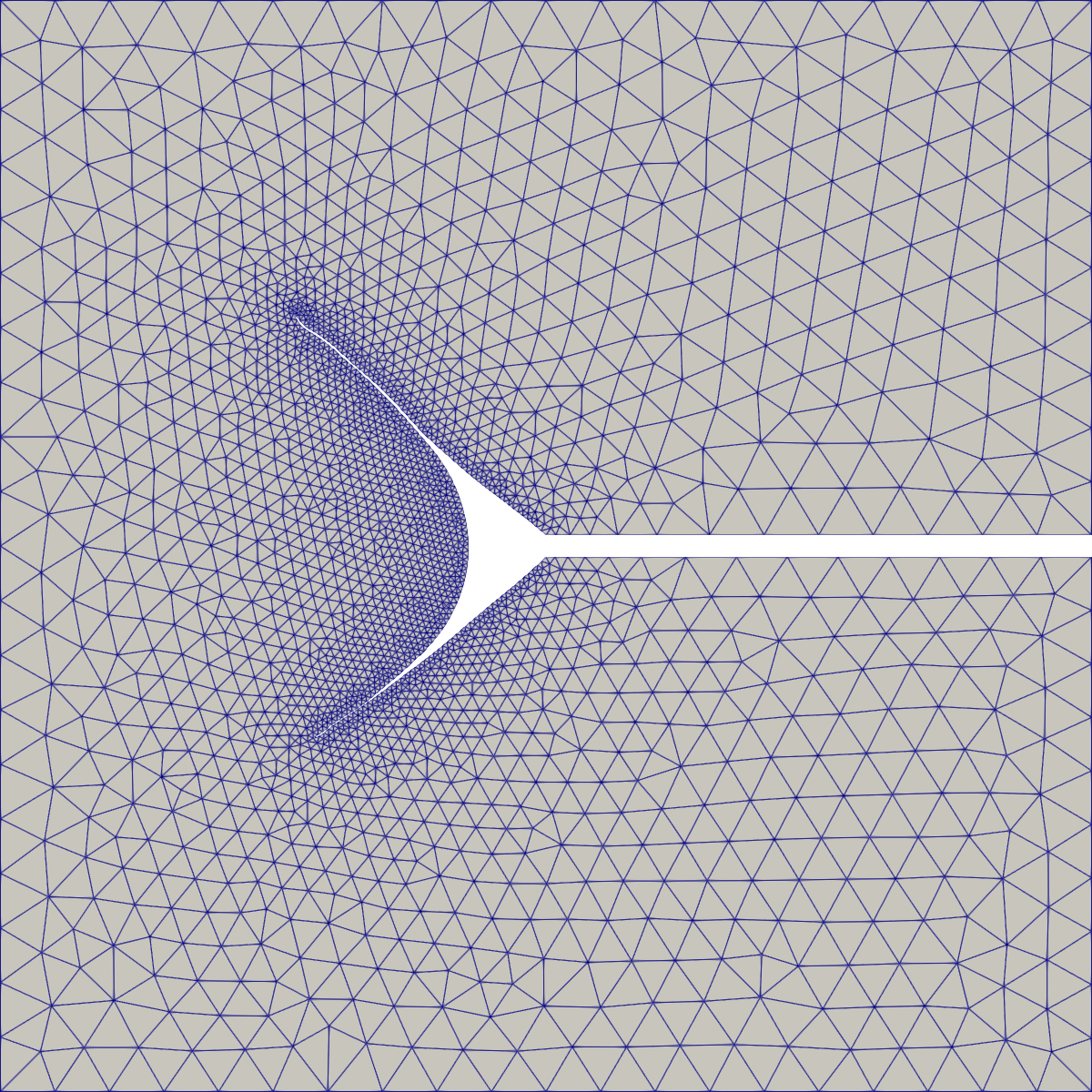}%
		\caption{The fracture continues to propagate with increasing displacement, shown here at $\SI{1.04e-2}{\milli\meter}$.}%
		\label{fig:mesh_shear_middle}%
	\end{subfigure}%
	\quad%
	\begin{subfigure}[t]{\figurewidth}%
		\centering%
		\includegraphics[width=\figurewidth]{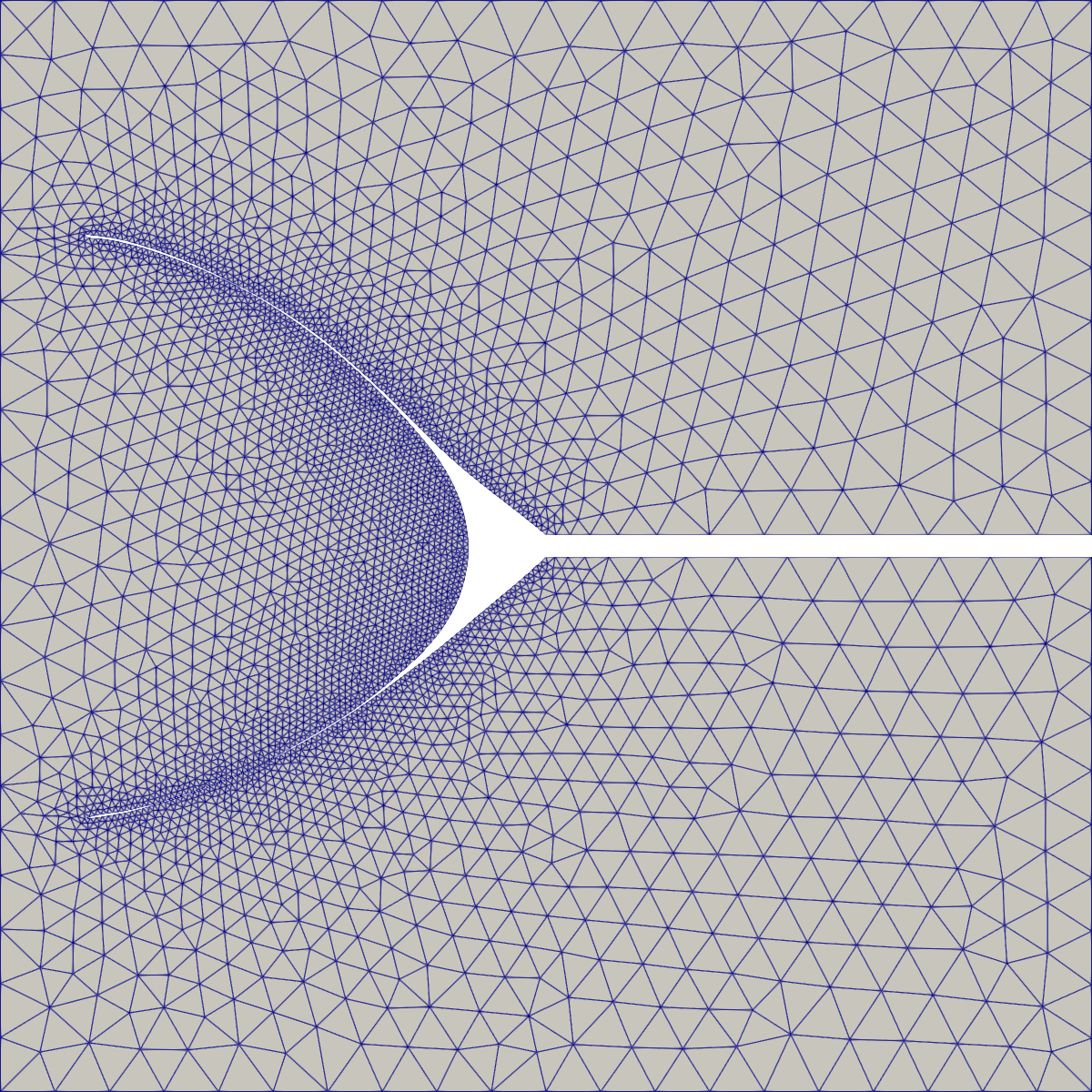}%
		\caption{Even at a displacement of $\SI{2.22e-2}{\milli\meter}$, the fracture has not completely propagated through the domain.}%
		\label{fig:mesh_shear_end}%
	\end{subfigure}%
	\caption{Fracture propagation for the shear test.}%
	\label{fig:mesh_shear}%
	\vspace*{\floatsep}%
	\centering%
	\setlength\figureheight{5cm}%
	\setlength\figurewidth{.45\textwidth}%
	\begin{subfigure}[t]{\figurewidth}%
		\centering%
		\includegraphics{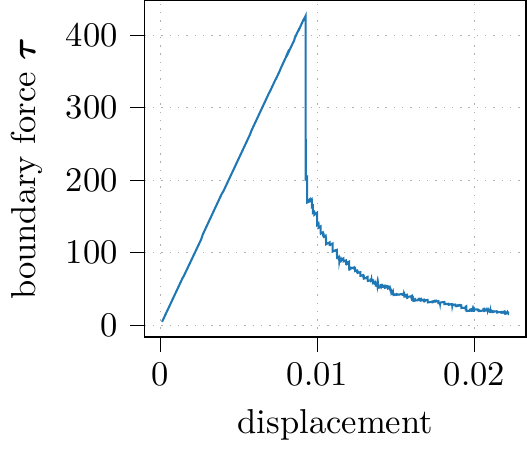}
		\label{fig:shear_boundary_force}%
	\end{subfigure}%
	\quad%
	\begin{subfigure}[t]{\figurewidth}%
		\centering%
		\includegraphics{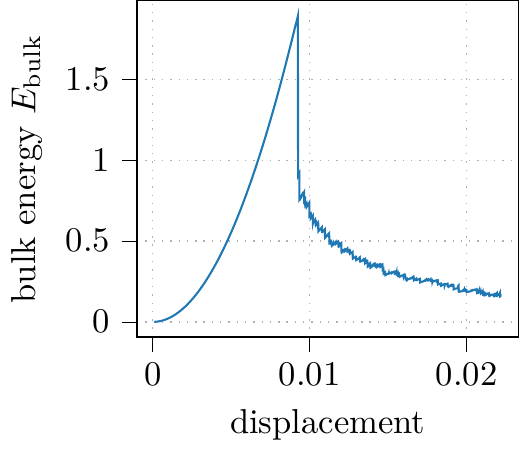}
		\label{fig:shear_bulkEnergy}%
	\end{subfigure}%
	\caption{Plot of boundary force $\vec{\tau}$ (left) and bulk energy $E_{\text{bulk}}$ (right) over displacement for the shear test.}%
	\label{fig:shear}%
	\vspace*{-2mm}
\end{figure}

The fracture propagation is depicted in Fig.~\ref{fig:mesh_shear}. The displacement where the first movement at the tip is detected occurs at $8.67\cdot 10^{-3}$, however the first major fracture propagation occurs at $9.26 \cdot 10^{-3}$.  The computational domain $\varOmega$ at the end of the gradient descent for a displacement of $9.26 \cdot 10^{-3}$ is shown in Fig.~\ref{fig:mesh_shear_initial}. The fracture is seen to propagate in the direction of the top-left and of the top-right corner simultaneously as no stress splitting is incorporated, see, e.g.,~\cite{AmbatiGerasimovLorenzis:2015}. With increasing displacement, the fracture propagates further, see Fig.~\ref{fig:mesh_shear_middle}, however it never completely splits the domain. The simulation is stopped at $2.22\cdot 10^{-2}$, for which the computational domain is shown in Fig.~\ref{fig:mesh_shear_end}. This mesh contains 3480 nodes and 6958 elements.

The initial linear increase in boundary force and quadratic increase in bulk energy is in accordance with literature, e.g.~\cite{AmbatiGerasimovLorenzis:2015,Mang2019}.
Similar to the tension test, we again find an initial major fracture growth at slightly lower displacement values than the isotropic results described in literature. Nonetheless, we also observe that the fracture gradually grows with increasing displacement. 

\section{Conclusion}
\label{sec:Conclusion}

In this paper, we have discussed a novel formulation for quasi-static fracture propagation by
means of shape optimization methods. Initial tests on two benchmark configurations yield
qualitative agreement with results obtained by phase-field formulations. Future studies need
to clarify the reasons for slight deviations in obtained load-displacements curves. Particular
questions would be the influence of the initially prescribed shape, i.e., the choice of
a rounded tip. Additional attention should be put on the choice of bilinear form for the deformation equation. Specifically, the choice of linear elasticity with parameters as in e.g.~\cite{SchulzSiebenbornWelker2015:2,Welker_diffeological} is demanding for the use in fracture propagation.

\paragraph*{Acknowledgements.}
This work has been partly supported 
by the state of Hamburg (Germany) within the Lan\-des\-for\-schungs\-för\-der\-ung under project\linebreak ``Simulation-Based Design Optimization of Dynamic Systems Under Uncertainties'' \mbox{(SENSUS)} with project number LFF-GK11.

\vspace{\baselineskip}

\bibliographystyle{hplain}
\bibliography{BibPaperCrackShapeOpt}

\begin{thebibliography}{10}

\bibitem{Alnaes2015}
M.~Alnæs, J.~Blechta, J.~Hake, A.~Johansson, B.~Kehlet, A.~Logg,
  C.~Richardson, J.~Ring, M.E. Rognes, and G.N. Wells.
\newblock The {FEniCS} {Project} {Version} 1.5.
\newblock {\em Archive of Numerical Software}, 3(100), 2015.

\bibitem{AmbatiGerasimovLorenzis:2015}
M.~Ambati, T.~Gerasimov, and L.~de~Lorenzis.
\newblock A review on phase-field models of brittle fracture and a new fast
  hybrid formulation.
\newblock {\em Comput. Mech.}, 55(2):383--405, 2015.

\bibitem{AmbrsosioTortorelli:1990}
L.~Ambrosio and V.M. Tortorelli.
\newblock Approximation of functionals depending on jumps by elliptic
  functional via {$\Gamma$}-convergence.
\newblock {\em Comm. Pure Appl. Math.}, 43(8):999--1036, 1990.

\bibitem{Artina:2015}
M.~Artina, M.~Fornasier, S.~Micheletti, and S.~Perotto.
\newblock Anisotropic mesh adaptation for crack detection in brittle materials.
\newblock {\em SIAM Journal on Scientific Computing}, 37(4):B633--B659, 2015.

\bibitem{BourdinFrancfortMarigot:2008}
B.~Bourdin, G.~A. Francfort, and J.-J. Marigo.
\newblock The variational approach to fracture.
\newblock {\em J. Elasticity}, 91(1-3):5--148, 2008.

\bibitem{BurkeOrtnerSuli:2010}
S.~Burke, C.~Ortner, and E.~S{\"u}li.
\newblock An adaptive finite element approximation of a variational model of
  brittle fracture.
\newblock {\em SIAM J. Numer. Anal.}, 48(3):980--1012, 2010.

\bibitem{Delfour-Zolesio-2001}
M.C. Delfour and J.-P. Zol\'esio.
\newblock {\em {Shapes and Geometries: Metrics, Analysis, Differential
  Calculus, and Optimization}}, volume~22 of {\em Adv. Des. Control}.
\newblock SIAM, 2nd edition, 2001.

\bibitem{FrancfortMarigo:1998}
G.A. Francfort and J.-J. Marigo.
\newblock Revisiting brittle fracture as an energy minimization problem.
\newblock {\em J. Mech. Phys. Solids}, 46(8):1319--1342, 1998.

\bibitem{GeiersbachHandbook}
C.~Geiersbach, E.~Loayza-Romero, and K.~Welker.
\newblock {PDE}-constrained shape optimization: {T}owards product shape spaces
  and stochastic models.
\newblock In: K.~Chen, C.-B. Sch\"{o}nlieb, X.-C. Tai, and L.~Younes, editors,
  {\em Handbook of Mathematical Models and Algorithms in Computer Vision and
  Imaging}. Springer, 2021.
\newblock Accepted for publication.

\bibitem{GeiersbachLoayzaWelker}
C.~Geiersbach, E.~Loayza-Romero, and K.~Welker.
\newblock Stochastic approximation for optimization in shape spaces.
\newblock {\em SIAM J. Optim.}, 31(1):348--376, 2021.

\bibitem{Geuzaine2009}
C.~Geuzaine and J.-F. Remacle.
\newblock Gmsh: A 3-d finite element mesh generator with built-in pre- and
  post-processing facilities.
\newblock {\em International Journal for Numerical Methods in Engineering},
  79(11):1309--1331, 2009.

\bibitem{Griffith:1921}
A.A. Griffith.
\newblock The phenomena of rupture and flow in solids.
\newblock {\em Philos. Trans. R. Soc. Lond.}, 221:163--198, 1921.

\bibitem{KrieglMichor}
A.~Kriegl and P.~Michor.
\newblock {\em The Convient Setting of Global Analysis}, volume~53 of {\em
  Mathematical Surveys and Monographs}.
\newblock American Mathematical Society, 1997.

\bibitem{LuftWelker}
D.~Luft and K.~Welker.
\newblock Computational investigations of an obstacle-type shape optimization
  problem in the space of smooth shapes.
\newblock In: {\em International Conference on Geometric Science of
  Information}, pages 579--588. Springer, 2019.

\bibitem{Mang2019}
K.~Mang, M.~Walloth, T.~Wick, and W.~Wollner.
\newblock Mesh adaptivity for quasi-static phase-field fractures based on a
  residual-type a posteriori error estimator.
\newblock {\em {GAMM}-Mitteilungen}, 43(1), 2019.

\bibitem{Michor2006}
P.W. Michor and D.~Mumford.
\newblock Riemannian geometries on spaces of plane curves.
\newblock {\em J. Eur. Math. Soc. (JEMS)}, 8:1--48, 2006.

\bibitem{miehe2010phase}
C.~Miehe, M.~Hofacker, and F.~Welschinger.
\newblock A phase field model for rate-independent crack propagation: {R}obust
  algorithmic implementation based on operator splits.
\newblock {\em Computer Methods in Applied Mechanics and Engineering},
  199(45-48):2765--2778, 2010.

\bibitem{MieheWelschingerHofacker:2010}
C.~Miehe, F.~Welschinger, and M.~Hofacker.
\newblock Thermodynamically consistent phase-field models of fracture:
  {V}ariational principles and multi-field {FE} implementations.
\newblock {\em Int. J. Numer. Meth. Eng.}, 83(10):1273--1311, 2010.

\bibitem{NeitzelWickWollner:2018}
I.~Neitzel, T.~Wick, and W.~Wollner.
\newblock An optimal control problem governed by a regularized phase-field
  fracture propagation model. part {II} the regularization limit.
\newblock {\em SIAM J. Control Optim.}, 3(57):1672--1690, 2019.

\bibitem{Schulz}
V.H. Schulz.
\newblock {A Riemannian view on shape optimization}.
\newblock {\em Found. Comput. Math.}, 14(3):483--501, 2014.

\bibitem{SchulzSiebenbornWelker2015:2}
V.H. Schulz, M.~Siebenborn, and K.~Welker.
\newblock Efficient {PDE} constrained shape optimization based on
  {S}teklov-{P}oincar{\'e} type metrics.
\newblock {\em SIAM J. Optim.}, 26(4):2800--2819, 2016.

\bibitem{SchulzWelker}
V.H. Schulz and K.~Welker.
\newblock On optimization transfer operators in shape spaces.
\newblock In: V.H. Schulz and D.~Seck, editors, {\em Shape Optimization,
  Homogenization and Optimal Control}, pages 259--275. Springer, 2018.

\bibitem{Siebenborn2017}
M.~Siebenborn and K.~Welker.
\newblock Algorithmic aspects of multigrid methods for optimization in shape
  spaces.
\newblock {\em SIAM J. Sci. Comput.}, 39(6):B1156--B1177, 2017.

\bibitem{Sokolowski1991}
J.~Sokolowski and J.~Zol\'esio.
\newblock {\em Introduction to Shape Optimization: Shape Sensitivity Analysis}.
\newblock Springer-Verlag, 1991.

\bibitem{WallothWollner:2021}
M.~Walloth and W.~Wollner.
\newblock A posteriori estimator for the adaptive solution of a quasi-static
  fracture phase-field model with irreversibility constraints.
\newblock {\em SIAM J. Sci. Comput.}, 44(3):B479--B505, 2022.

\bibitem{Welker2016}
K.~Welker.
\newblock {\em Efficient {PDE} Constrained Shape Optimization in Shape Spaces}.
\newblock PhD thesis, Universit\"{a}t Trier, 2016.

\bibitem{Welker_diffeological}
K.~Welker.
\newblock Suitable spaces for shape optimization.
\newblock {\em Applied Mathematics and Optimization}, 2021.
\newblock DOI: 10.1007/s00245-021-09788-2.

\bibitem{WheelerWickWollner:2014}
M.F. Wheeler, T.~Wick, and W.~Wollner.
\newblock An augmented-{L}agrangian method for the phase-field approach for
  pressurized fractures.
\newblock {\em Comput. Methods Appl. Mech. Engrg.}, 271(1):69--85, 2014.

\end{thebibliography}

\end{document}